 \documentclass[preprint,review,12pt]{elsarticle}

\usepackage[latin1]{inputenc}
\usepackage[english]{babel}
\usepackage{bm}
 \usepackage{epsfig}

\usepackage{amssymb}

\usepackage[nodots]{numcompress}

\newcommand{\ds}{\displaystyle}

\def\div{\hbox{div\,}}

\newcommand{\bea}{\begin{eqnarray}}
\newcommand{\eea}{\end{eqnarray}}
\newcommand{\nn}{\nonumber}
\newcommand{\ddt}{\frac{\partial}{\partial t}}
\newcommand{\ddx}{\frac{\partial}{\partial x}}
\newcommand{\V}{{\bm V}}
\newcommand{\F}{{\bm F}}

\newcommand{\Vol}{\textrm{Vol}}
\newcommand{\dVol}{\delta \textrm{Vol}}

\journal{The European Journal of Mechanics B Fluids}

\begin{document}

\begin{frontmatter}



\title{A totally Eulerian Finite Volume solver for multi-material
fluid flows : Enhanced Natural Interface Positioning (ENIP)}


\author[Raph]{Raphaël Loubère\corref{cor1}}
\ead{raphael.loubere@math.univ-toulouse.fr}
\address[Raph]{CNRS et Universit\'e de Toulouse
 IMT (Institut de Math\'ematiques de Toulouse),
 31062 Toulouse, France}

\author[Pilou]{Jean-Philippe Braeunig}
\ead{braeunig@math.u-strasbg.fr}
\address[Pilou]{INRIA Nancy Grand-Est, Equipe CALVI,  615 rue du Jardin Botanique 54600 Villers-l\`es-Nancy, France\\
CEA DIF Bruy\`eres-le Ch\^atel, 91297 Arpajon, France}

\author[JMG]{Jean-Michel Ghidaglia}
\ead{jmg@cmla.ens-cachan.fr}
\address[JMG]{CMLA, CNRS et ENS de Cachan
61 Av. du Pr{\'e}sident Wilson
Cachan Cedex 94235, France}

\cortext[cor1]{Corresponding author}


\begin{abstract}

This work concerns the simulation of compressible multi-material fluid flows and follows the method FVCF-NIP described in the former paper \cite{Braeunig09}. This Cell-centered Finite Volume method is totally Eulerian since the mesh is not moving and a sharp interface, separating two materials, evolves through the grid. A sliding boundary condition is enforced at the interface and mass, momentum and total energy are conserved. Although this former method performs well on 1D test cases, the interface reconstruction suffers of poor accuracy in conserving shapes for instance in linear advection. This situation leads to spurious instabilities of the interface. The method Enhanced-NIP presented in the present paper cures an inconsistency in the former NIP method that improves strikingly the results. It takes advantage of a more consistent description of the interface in the numerical scheme.  Results for linear advection and compressible Euler equations for inviscid fluids are presented to assess the benefits of this new method.
\end{abstract}

\end{frontmatter}


\tableofcontents

\section{Introduction} \label{sec:introduction}

 The two-material compressible hydrodynamics equations (Euler equations) are considered
in this work.
The flow regime is such that molecular viscosity within materials is
neglected: materials are supposed immiscible and separated by sharp
interfaces, with perfect sliding between materials.
Each material is characterized by its own equation of state (EOS).

  The formalism of finite volume methods is close to the mechanical
viewpoint, and generic for different types of physical models.
  Thus, it might be easier to add such models; surface tension
or turbulent diffusion for instance.
  The discretization order is limited, but this method is accurate
to simulate hydrodynamic shock waves, because of the consistency between numerical
treatment and mechanics.\\
  The extension of Eulerian schemes to multi-material fluid flows can be
obtained by various techniques.
  One is to introduce the cell mass fraction $c_\alpha$ of material
$\alpha$ and let it evolve according to material velocity.
  The cell is called pure if a material $\alpha$ satisfies
$c_\alpha=1$ and is called mixed if $c_\alpha \in ]0,1[$.
  Pure cells filled by material $\alpha$ are calculated in the same
manner as for the single material method.
  Mixed cell evolution is computed using a mixing equation
of state that takes into account material mass fractions, see {\it e.g.} \cite{int_diff1}.
One drawback of this approach is the numerical diffusion
of the interface that obviates sharp interface capturing. It turns out that for some applications, this drawback is not acceptable since the diffusion of one material into another one will correspond to a different physics. For example the two material could react when a molecular mixture is formed. Hence such a diffusion should occur only for physical reasons and not for numerical ones.

 In the case of sharp interface capturing methods, the interface is approximated
in a mixed cell by a segment by most authors.
However more complex curves than straight line or more complex theory (see \cite{MOF} for instance)
might be used.
  A famous method using sharp interface reconstruction is the
Lagrange+Remap Finite Volume scheme, initiated in \cite{NohWood} and
further improved in \cite{DLY1}. It belongs to the family of so called
Volume of Fluid (VOF) methods.
  The first step of this method is a Lagrangian scheme, resulting in a mesh displacement with
 material velocity.
  The second step is a multi-material remapping of Lagrangian mesh onto the
original Eulerian mesh, by exchanging volume fluxes between cells
related to the Lagrangian motion of cell edges.
  The new interface position in mixed cells is determined using
the partial volumes of the materials and the interface normal vector.
The later is
calculated using volume fractions from neighboring cells.
  Thus the ratio of each material in volume fluxes is
deduced from the multi-material remapping.
  Some methods with the same kind of operator splitting are used
for incompressible multi-material fluid flows as in
\cite{Zaleski}.
  These methods provide sharp interface between materials and discontinuous
quantities in mixed cells, allowing large deformations and
transient flows.
  In this context, the drawback of these Lagrange+Remap methods is the limited accuracy of
the underlying single phase scheme due to diffusion
induced by the remapping step.
Moreover, more complex physics at material interfaces such as
sliding effects, is not possible.

The FVCF scheme (Finite Volume with Characteristic Flux) has been introduced in \cite{GKL1} for simulating single phase compressible flows or multi-phase models without sharp interface capturing.
 The method described in \cite{Braeunig09}, so called NIP method (Natural Interface Positioning),
is an add-on to the FVCF method in order to deal with multi-material fluid flows with sharp interface capturing.
It is a cell centered totally Eulerian scheme, in which
material interfaces are represented by a discontinuous piecewise
linear curve. A treatment for interface evolution is proposed on Cartesian
structured meshes which is locally conservative in mass, momentum and total energy and
allow the materials to slide on each others.
Discrete conservation laws are written on
partial volumes as well as on pure cells, considering the interface in the cell as a
moving boundary without any diffusion between materials.
A specific data structure called
{\em condensate} is introduced in order to write a finite volume scheme even when the
considered volume is made of moving boundaries, i.e. interfaces. This treatment includes
an explicit computation of pressure and velocity at interfaces.

In \cite{Braeunig09} are shown $2D$
results illustrating the capability of the method to deal with perfect sliding, high pressure
ratios and high density ratios.
This former method however produces non satisfactory results in the context of advection of geometrical shapes
especially when dealing with low Mach numbers.
It is however a classical misbehavior of most of advection and reconstruction methods  which have a tendency
to destroy the shape of advected objects due to numerical approximations.
However, this former method gives very poor results when advecting geometrical shapes especially when dealing with low Mach number flows.
In this work we propose a new method called ENIP (Enhanced NIP) that is an improvement of  the NIP method by a more accurate treatment of condensates. On a very simple
example: the advection of a square, an inconsistency in the NIP interface reconstruction method will be exhibited. We will then introduce ENIP that cures this situation.
Numerical examples are presented in the last Section to assess the validity and efficiency
of this new approach.

\section{FVCF-ENIP: Finite Volume Characteristics Flux with
Enhanced Natural Interface Positioning technique}
\label{sec:VFFC-NIP}

\subsection{Governing equations}
The model addressed  in this work is the compressible Euler equations in space dimension $d$ that can be
written in a conservative form as follows:
\bea
\label{euler_equation_rho}
 \ddt (\rho) + \div (\rho u)             &=& 0, \\
\label{euler_equation_rhou}
\ddt (\rho u)
  + \div (\rho u \otimes u + p I) &=& 0, \\
\label{euler_equation_rhoE}
 \ddt (\rho E)
  + \div ((\rho E+p)u )                &=& 0,
\eea
where $\rho$ denotes the density, $u \in \mathbb{R}^d$ the velocity field,
 $p$ the pressure, $E=e+| u |^2/2$
the specific total energy and $e$ the specific internal energy.
  An equation of state of the form $EOS(\rho,e,p)=0$ or $p=p(\rho, e)$
is provided in order to close the system.\\
Let us consider a generic
conservative form with ${\V}=(\rho,\rho u,\rho E)^t$ the unknown
vector of conservative variables and flux $\F$ is a matrix
valued function defined as:
\begin{equation}
\begin{array}{llcl}
\F :  & \mathbb{R}^{d+2} &  \longrightarrow &  \mathbb{R}^{d+2} \times \mathbb{R}^{d}\\
      & \V               &  \longmapsto     & \F(\V),
\end{array}
\end{equation}
for all direction $n \in \mathbb{R}^{d}$, $\F(\V) \cdot n$ is given
in terms of $\V$ by:
\begin{equation}
    \F(\V) \cdot n= \left( \rho  \left(u \cdot n\right),
       \rho u  \left(u \cdot n\right) + p n,
        \left(\rho E +p\right)  \left(u \cdot n\right) \right).
\end{equation}
The compressible Euler equations  (\ref{euler_equation_rho}-\ref{euler_equation_rhoE}) can then
be rewritten as:
\begin{equation}
    \partial_t \V + \div \F(\V)  = 0.
\end{equation}

\subsection{FVCF: Single material scheme}
FVCF method uses a directional splitting on Cartesian
structured meshes. The method is thus detailed for only one generic direction denoted by $x$.
In $d$
dimensions of space, the algorithm described for direction $x$ has
to be replicated $d$ times, one for each direction. However, this directional splitting does
not modify at all the underlying single material scheme FVCF for pure cells.
 In \textit{2D}:
\begin{itemize}
  \item[-] variables at $t^{n,x}$ are
calculated from those at $t^n$ by the $x$ direction  step,
  \item[-] variables at $t^{n+1}$ are calculated
from those at $t^{n,x}$ by the $y$ direction  step.
\end{itemize}
\bea
\label{2step1}
  \Vol_i \frac{\V_i^{n,x }-\V_i^{n }}{\Delta t} + A_x \left(\bm{\phi}^n_\ell+\bm{\phi}^n_r \right) &=&  0,\\
\label{2step2}
  \Vol_i \frac{\V_i^{n+1}-\V_i^{n,x}}{\Delta t} + A_y \left(\bm{\phi}^n_d   +\bm{\phi}^n_u \right) &=&  0,
\eea
where the cell volume is $\Vol_i$, the cell face area are $A_x$ and $A_y$ respectively
normal to $x$ and $y$ directions,
up, down, right and left direction fluxes $\bm{\phi}^n_u$, $\bm{\phi}^n_d$, $\bm{\phi}^n_r$, $\bm{\phi}^n_\ell$
calculated with respect of the outgoing normal direction $n_d$ of cell face
$\Gamma_d$ in direction $d$ using variables at time $t^n$, i.e.
\bea
\bm{\phi}^n_d=\frac{1}{A_d} \int_{\Gamma_d} \F(\V^n) \cdot n_d dS.
\eea
This flux is further approximated using the finite volume scheme FVCF described in \cite{GKL1}.

\subsection{FVCF-NIP: Multi-material scheme}

One considers multi-material flows. The subcell model addressed here for the multi-material representation is
a cell $C$ of volume $\Vol_C$ containing $n_m$ different materials, each of them filling
a partial volume $\Vol^k_{C}$ such that
\bea
 \ds \sum^{n_m}_{k=1} \Vol^k_{C}=\Vol_{C}.
\eea
 Cell $C$ is referred to as pure if $n_m=1$,
and as mixed if $n_m>1$. The interfaces in mixed cells are  approximated by segments separating materials into two partial volumes which are pure on both sides of the interface. \\
A partial volume cell-centered variable vector
$\V_{k}=(\rho_k,\rho_k u_k,\rho_k E_k)^t$ and an equation of state
$EOS_{k} (\rho_k,e_k,p_k)=0$ are also associated with each
material labeled by $k \leq n_m$ in the mixed cell.\\
 FVCF-NIP method uses a directional
splitting scheme for the interface evolution without loosing the accuracy of the
Eulerian scheme in the bulk of materials.
Consequently this scheme is restricted to structured Cartesian
mesh.

\begin{figure}
  \begin{center}
  \epsfxsize=0.485\textwidth
  \epsfysize=0.27\textwidth
  \epsfbox{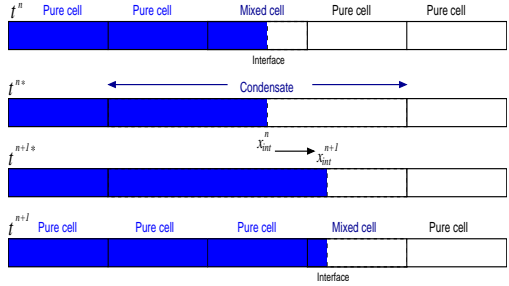}
  \caption{ \label{fig_evol_cond}
    Sketch of a condensate. Evolution of an interface through a cell edge
    during one time step. Mixed and pure neighbor cells are merged to obtain
    the so called condensate at fictitious time $t^{n*}$. Interface evolution is performed
    within this condensate from $t^{n*}$ to $t^{n+1*}$
    This condensate is then split back into Eulerian cells.
  }
  \end{center}
\end{figure}

The multi-material extension proposed in \cite{Braeunig09} considers the finite
volume scheme (\ref{2step1}-\ref{2step2}) on each partial volume in a mixed
cell. The obtained scheme is conservative by construction and is constrained
with the same \textit{CFL} condition as the single material
scheme\footnote{Without such a special treatment the time step would be constrained
by the smallest partial volume, which is arbitrarily small.}.
  NIP method consists in removing cell edges when this cell contains
an interface.
  Therefore each partial volume is merged with the
neighbor pure cells filled with the same material, see
Figure \ref{fig_evol_cond}.
Variables in these enlarged partial
volumes are obtained by writing the conservation laws on the merged volumes
\bea
\overline{\Vol_1} &=& \Vol_1+\Vol_{pure~1}, \\
\overline{\Vol_2} &=& \Vol_2+\Vol_{pure~2},
\eea
then on the conserved variables
\bea
\overline{\V}_1  &=& \frac{\Vol_1 ~\V_1 + \Vol_{pure~1}~ \V_{pure~1}}{\overline{\Vol_1}} \\
\overline{\V}_2  &=& \frac{\Vol_2 ~\V_2 + \Vol_{pure~2}~ \V_{pure~2}}{\overline{\Vol_2}}.
\eea
This set of cells is associated with its left and
right single material fluxes $\bm{\phi}_\ell$ and $\bm{\phi}_r$.
  Internal cell edges are forgotten, considering only
enlarged volumes $\overline{\Vol_1}$ and $\overline{\Vol_2}$ and
averaged variables $\overline{\V}_1$ and $\overline{\V}_2$,
separated by an interface; this system is called a condensate.

  Actually, this numerical strategy consists in condense
neighboring mixed cells in one direction of the Cartesian mesh, in which
interfaces are considered as mono dimensional objects, namely they are
considered vertical during $x$ direction step and horizontal during $y$ direction step.
A condensate then contains layers of successive different
materials that are separated by straight interfaces. The thickness of
these layers is calculated through volume conservation. The ordering of
layers is given by the \textit{2D} description from the
previous time step. It is determined thanks to the volume fractions of
neighboring cells.
  The layer evolution is calculated in a Lagrangian fashion
which implies that layers can be as thin as partial volumes are small.
Once quantities and interface positions inside the
condensate are known at time $t^{n+1}$, they are remapped back
onto the original Eulerian mesh. Finally a 2D normal in each mixed cell
is computed as described in \cite{DLY1}: the method is based on an approximation
of the gradient of the volume fraction function
in mixed cells.
It provides the normal to materials interface in each cell that is further used to locate materials
within mixed cells.
The numerical scheme used in a condensate
is presented in great details in \cite{Braeunig09}
and we omit this description in this work and rather focus on the interface reconstruction method.

%
As shown in \cite{Braeunig09} this numerical method has several attractive properties as
conservation and perfect sliding of materials as instance.
Moreover $\Delta t$ is not restricted by small partial volume thanks to a tight control of density
and pressure \cite{Braeunig10}.
The numerical experiments carried out in \cite{GKL1,these_jpb,Braeunig09,benchmark}
have confirmed the efficiency of such a method for compressible multi-material computation.
Although very promising, the method suffers from the way interfaces are dealt with.

\begin{figure}
\begin{center}
	\epsfxsize=0.485\textwidth
        \epsfbox{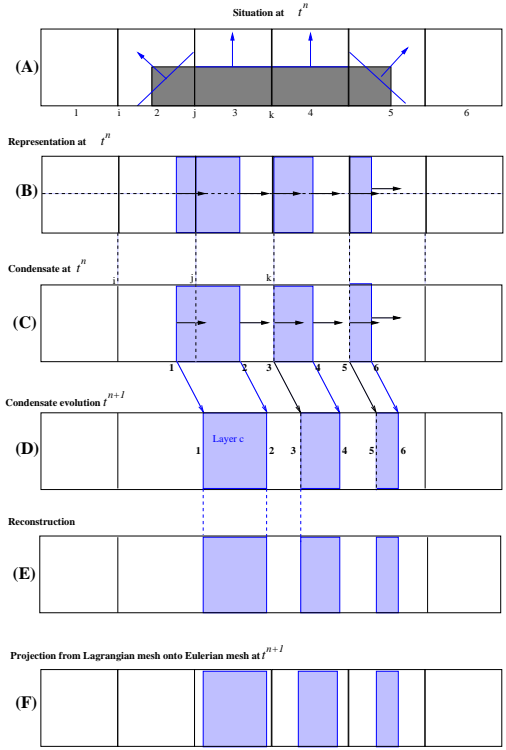}
\caption{ \label{fig:nip}
NIP method ---
\textbf{(A)} Situation at $t^n$ with real materials geometry, interfaces and normals to them.
\textbf{(B)} Representation of partial volumes at $t^n$.
\textbf{(C)} Construction of a condensate at $t^n$ by merging layers of contiguous partial volumes of the same material.
\textbf{(D)} Evolution of condensate in a Lagrangian fashion during $\Delta t$.
\textbf{(E)} Condensate reconstruction at $t^{n+1}$.
\textbf{(F)} Condensate projection/remapping from Lagrangian mesh onto original mesh.}
\end{center}
\end{figure}

\begin{figure}
\begin{center}
\hspace{-1cm}
	\epsfxsize=0.28\textwidth
        \epsfbox{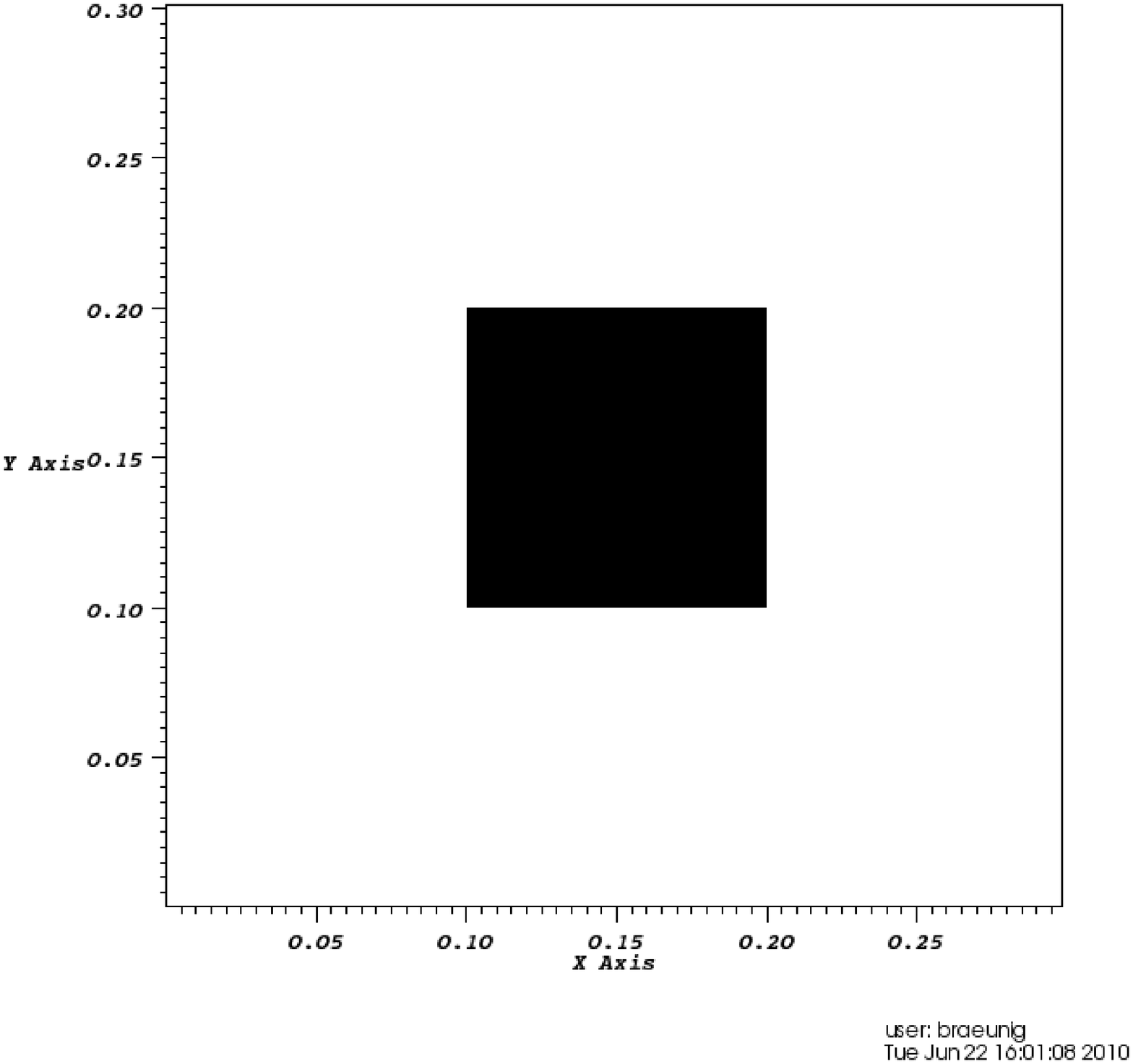}
\hspace{-1cm}
	\epsfxsize=0.28\textwidth
        \epsfbox{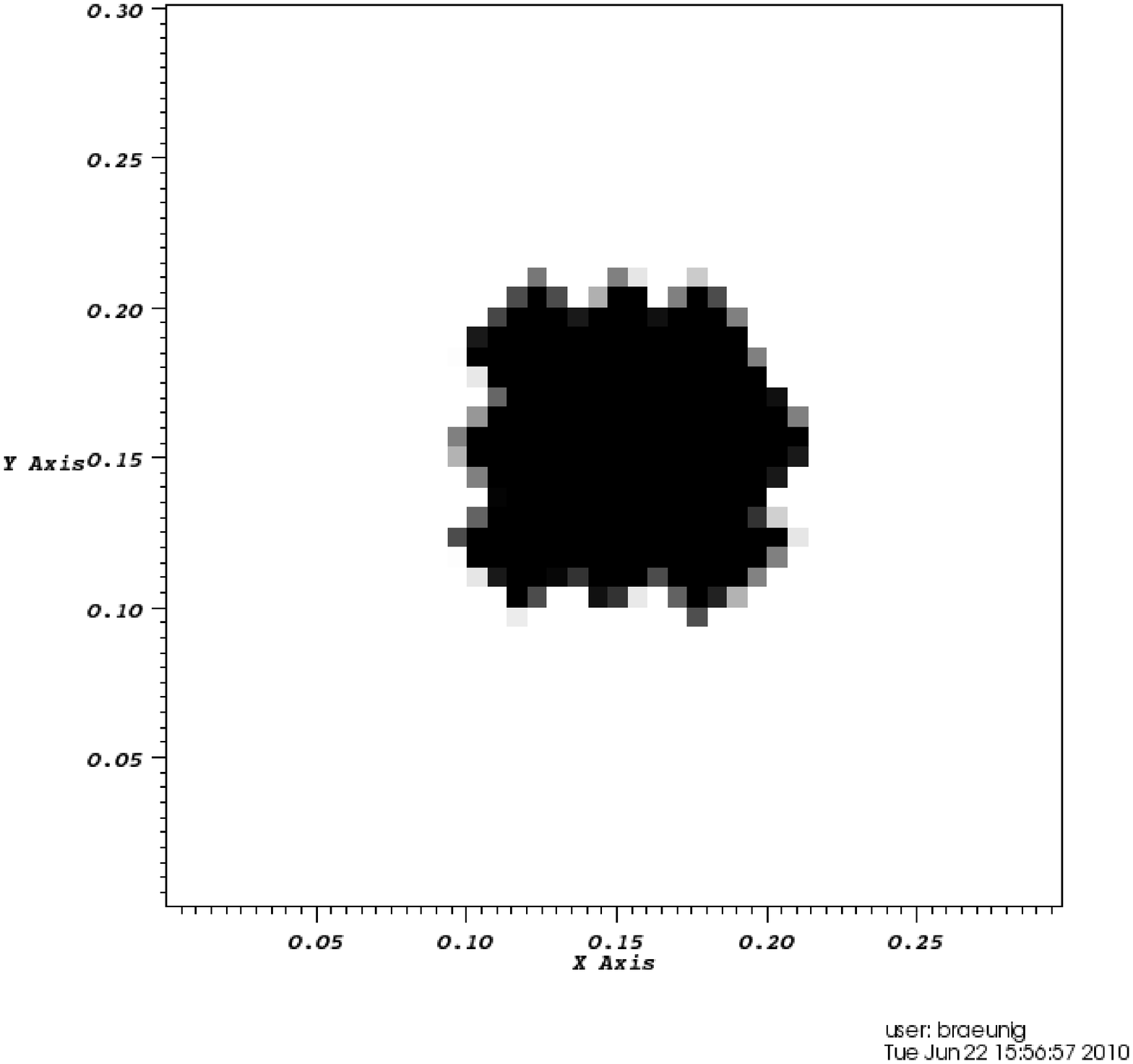}
\caption{ \label{fig:advection_square}
Diagonal advection of a square by NIP method ---
Left: Exact configuration ---
Right: Numerical configuration obtained by NIP. }
\end{center}
\end{figure}

In order to illustrate the interface reconstruction method
NIP let us consider a square like interface cutting the Eulerian cells,
as in Figure \ref{fig:nip}-(A). These interfaces are indeed defined by their normals within each cell.
NIP method consists of the following steps assuming the condensate is in the $x$ direction:
\begin{itemize}
\item {\em Representation} Figure \ref{fig:nip}-(B).
The representation step can be seen as the way of determining on which side (left or right)
of the mixed cell the material is to be put. This is done by comparing the direction of the
interface normal at time $t^n$ with the vertical direction.
\item {\em Condensate construction} Figure~\ref{fig:nip}-(C).
The construction of the condensate consists in discarding any cell edges in the mixed cells considered.
Then the partial volumes of the same contiguous materials are glued together into so called condensate layers.
As instance cell $2$ and $3$ dark materials are merged into one stand-alone layer with associated volume averaged values.
\item {\em Condensate evolution} Figure~\ref{fig:nip}-(D).
The condensate layers evolution is computed from $t^n$ to $t^{n+1}$
thanks to the numerical scheme developed in \cite{Braeunig09}. In short, each vertical interface
is assigned a velocity and, consequently, a new position of each layer within the condensate
is determined in a Lagrangian way. Any conserved variable is computed accordingly.
\item {\em Reconstruction} Figure~\ref{fig:nip}-(E).
This phase consists in ``guessing'' the shape of each layer in the condensate
before remapping.
The reconstruction phase was not originally considered as a true phase of the algorithm
as the author used the same shapes as the ones produced in phase {\em Condensate construction}, i.e. only vertical interfaces.
\item {\em Projection} Figure~\ref{fig:nip}-(F).
The projection step consists in remapping the shapes obtained from the reconstruction phase onto
the Eulerian grid. This step produces updated partial volumes in mixed cells. Volume
fractions are deduced.
\end{itemize}
When all mixed cells in the domain are treated for direction $x$,  the interface normals are computed using the updated volume fractions. This concludes the system evolution in direction $x$, as we are back to a similar situation as the one described in Figure~\ref{fig:nip}-(A).\\
In the case where the normal is almost vertical, positioning the material
on either side of the cell might be, at least inaccurate, or, worse, incorrect.
Furthermore the reconstruction phase is here clearly
inconsistent: the interfaces are initially horizontal in cell 3 and 4  (Figure~\ref{fig:nip}-(A)), while in  the Reconstruction Figure~\ref{fig:nip}-(E) and in the Projection Figure~\ref{fig:nip}-(F) phase interfaces are set vertical for any initial geometry.  This situation of a horizontal interface is the worst case, but it illustrates the lack of geometrical consistency of NIP. This inaccurate reconstruction
step leads to a lack of accuracy of the volume fractions obtained after the remapping step.
Ultimately, it impacts the whole numerical method in any advection process. \\
As an illustration let us consider the diagonal advection of a square back and forth as shown in
Figure~\ref{fig:advection_square}. We omit the exhaustive description of this test as it will be done
in the numerical Section of this paper. On the right panel it is obvious that the shape of the square is not
well approximated. More important the horizontal and vertical edges of the square
do not remain so.  This behaviour is less pronounced if one refines the mesh but still remains.
Our goal is to improve the reconstruction step so that the new method, denoted ENIP standing for Enhanced Natural Interface Positioning, cure this geometrical inconsistency during the advection phase of the algorithm.

\subsection{FVCF-ENIP}

\begin{figure}
\begin{center}
	\epsfxsize=0.485\textwidth
        \epsfbox{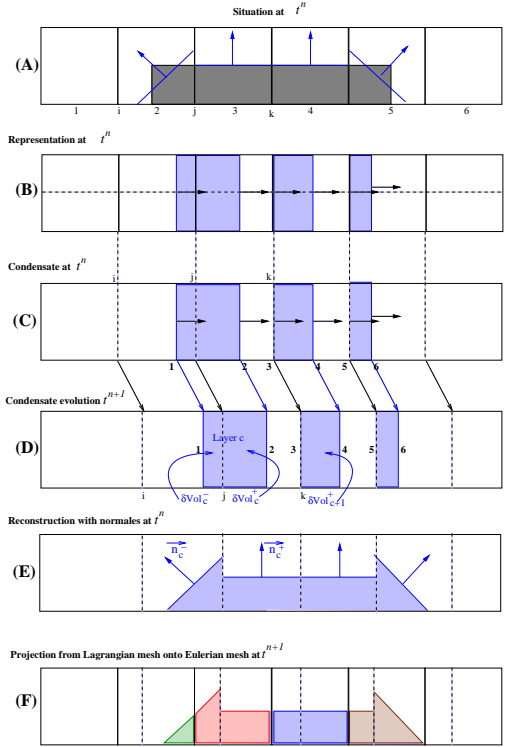}
\caption{ \label{fig:nip-2}
ENIP method ---
\textbf{(A)} Situation at $t^n$ with real material, interfaces and normals to them.
\textbf{(B)} Representation of material at $t^n$.
\textbf{(C)} Construction of a condensate at $t^n$ by merging of mixed cells leading
to layers of contiguous pieces of the same material.
\textbf{(D)} Evolution of condensate in a Lagrangian fashion during $\Delta t$.
Determine layer compression rates $\dVol_c^\pm$ through the evolution
of Lagrangian cells during $\Delta t$.
\textbf{(E)} Condensate reconstruction at $t^{n+1}$ using interface normals defined at $t^n$.
\textbf{(F)} Condensate projection/remapping from Lagrangian mesh onto Eulerian mesh.}
\end{center}
\end{figure}

The main idea of the new interface reconstruction method ENIP emanates from the following remarks:
\begin{enumerate}
\item At time $t^n$ any interface normal in mixed cell $i$ denoted $\vec{n}_i$ is known. It is
used to locate the partial volumes within cell $i$ when the condensate is constructed (phase
(B) and (C) of Figure~\ref{fig:nip}). However $\vec{n}_i$ is never taken into account in the reconstruction
and projection phases (E) and (F) from the same figure.
\item  Any layer of the condensate evolves as a Lagrangian object in the original method. Consequently the cell faces could evolve in an almost Lagrangian manner within this condensate. This makes possible to conserve the initial geometry of partial volumes during this Lagrangian motion.
\end{enumerate}
Therefore ENIP modifies several steps of NIP as depicted in Figure~\ref{fig:nip-2}.
Once a patch of neighbor mixed cells in $x$ direction\footnote{The  $y$ direction is treated likewise.} are agglomerated,
The same five steps as for NIP method are performed.
The first two steps are kept unmodified. The last three are modified as described
in the following.

\subsubsection{Lagrangian {\em Condensate evolution} step}
\begin{figure}
\begin{center}
	\epsfxsize=0.489\textwidth
        \epsfbox{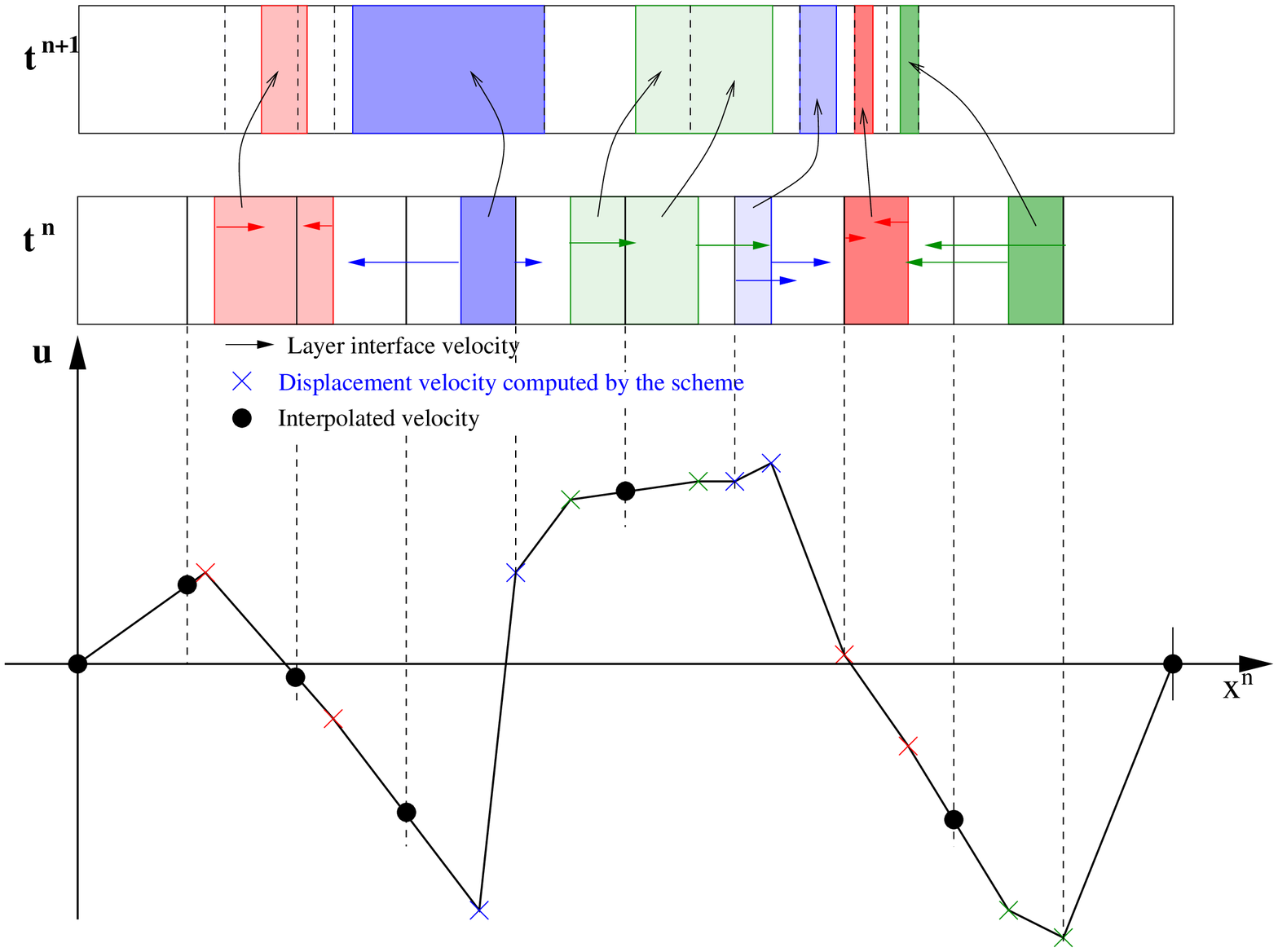}
\caption{ \label{fig:vitesse}
Sketch of linear displacement assumption ---
Displacement velocity varies linearly between the layer interface velocities
($\times$ in color) computed by the numerical scheme.
The cell interface velocity ($\bullet$) is interpolated.
The two top rows represent the evolution of a condensate in the $x$ direction
from $t^n$ to $t^{n+1}$. }
\end{center}
\end{figure}

\paragraph{Cell interface Lagrangian velocity}
After the condensate at $t^n$ is constructed,
each layer labeled $c$ is located thanks to the
left and right interface position respectively called $x_c^-, x_c^+$. The numerical scheme
provides the layer evolution, and as a by-product, the velocity of these interface positions,
$u_c^-,u_c^+$ are given by
\bea
x_c^{-,n+1} = x_c^-   + \Delta t \ u_c^- ,\ \ \ \
x_c^{+,n+1} = x_c^+   + \Delta t \ u_c^+ .
\eea
We make the following fundamental linear displacement assumption:
\textit{The velocity linearly varies within any layer}, see Figure~\ref{fig:vitesse}
for a sketch.
This assumption implies that any point $x_i \in \left[x_c^-;x_c^+\right]$
characterized by its 1D barycentric coordinates
\bea
\lambda_i^- = \frac {x_c^+-x_i}{x_c^+-x_c^-},  \  \ \ \ \   \lambda_i^+ = \frac {x_i-x_c^-}{x_c^+-x_c^-},
 \eea
moves to location
\bea\label{vit_noeud}
x_i^{n+1} = \lambda_i^- x_c^{-,n+1} +  \lambda_i^+ x_c^{+,n+1}
          = x_i + \Delta t \left( \lambda_i^- u_c^{-} +  \lambda_i^+ u_c^{+} \right) .
\eea
Then the point velocity is naturally set to
$u_i= \lambda_i^- u_c^{-} +  \lambda_i^+ u_c^{+}$.
Using this previous formula one can associate a ``Lagrangian''
velocity to any cell interface.
As instance in Figure~\ref{fig:nip-2}-(C) cell interface located at $x_i^n$ moves
to the position $x_i^{n+1}=x_i^n + \Delta t \ u_i$ with
$u_i$ being the linear combination between $u_c^{-}$ and $u_c^{+}$ {\it via} the
barycentric coordinates of point $x_i$ in $[x_c^-;x_c^+]$.
With the same formula one gets $x_{i+1}^{n+1} = x_{i+1}^n + \Delta \ u_{c+1}^-$
in the next layer as $x_{i+1}^n \equiv x_{c+1}^-$.
\paragraph{Compression/expansion rates}
The global rate of compression/expansion in layer $c$ during $\Delta t$ is given by
\bea
\dVol_c = \frac{x_c^{+,n+1}-x_c^{-,n+1}}{x_c^{+}-x_c^{-}}
           = 1 + \Delta t \frac{u_c^{+}-u_c^{-}}{x_c^{+}-x_c^{-}}.
\eea
The linearity assumption provides a simple way to determine
the rates of compression/expansion at left/right of a point
$x_i\in \left[x_c^-;x_c^+\right]$
\bea
\dVol_c^- = \frac{x_i^{n+1}-x_c^{-,n+1}}{x_c^{+}-x_c^{-}},\ \ \
\dVol_c^+ = \frac{x_c^{+,n+1}-x_i^{n+1}}{x_c^{+}-x_c^{-}},
\eea
that fulfil $\dVol_c^- + \dVol_c^+  = \dVol_c$.
Moreover the substitution of $x_i^{n+1}$ in the previous equations
yields
\bea
\dVol_c^- = \frac{x_i-x_c^-}{x_c^+ -x_c-} +
               \Delta t \frac{u_i-u_c^{-}}{x_c^{+}-x_c^{-}}
             =  \lambda_i^++ \Delta t \frac{u_i-u_c^{-}}{x_c^{+}-x_c^{-}},
\eea
where $u_i-u_c^- = (\lambda_i^- u_c^{-} +  \lambda_i^+ u_c^{+}) - u_c^-
= \lambda^+_i (u_c^+-u_c^{-})$, therefore the compression rates simply
writes
\bea
 \dVol_c^- &=&  \lambda_i^+ \left( 1 + \Delta t \frac{u_c^+-u_c^{-}}{x_c^{+}-x_c^{-}} \right)
               = \lambda_i^+ \dVol_c ,\\
 \dVol_c^+ &=& \lambda_i^-  \left( 1 + \Delta t \frac{u_c^+-u_c^{-}}{x_c^{+}-x_c^{-}} \right)
               =\lambda_i^- \dVol_c .
\eea
Each $\dVol_c^+$ or $\dVol_c^-$ is associated to a unique Eulerian cell; as instance
in Figure~\ref{fig:nip-2},  $\dVol_c^-$ is associated to cell $2$, $\dVol_c^+$ to
cell $3$, $\dVol_{c+1}^+$ to cell $4$ and so on.
Therefore $\dVol_c^\pm$ provides {\it de facto}
the compression/expansion of the partial volume originating from its associated Eulerian
cell motion. Furthermore, as any Eulerian mixed cell $i$ possesses a unique normal denoted
$\vec{n}_i$, this last is associated to the corresponding partial volume  $\dVol_c^\pm$;
this normal is consequently labeled $\vec{n}_c^\pm$.
These rates are then used to reconstruct the material topology into the
Lagrangian cell.

\subsubsection{{\em Reconstruction} step}
The Lagrangian cell $i+1/2$ at $t^{n+1}$ the interfaces of which moved as
\bea
x_i^{n+1} = x_i + \Delta t \ u_i , \ \ \ \ x_{i+1}^{n+1} = x_{i+1} + \Delta t \ u_{i+1},
\eea
changed its volume as
\bea
\dVol_{i+1/2} =  \frac{V_{i+1/2}^{n+1}}{V_{i+1/2}} = \frac{x_{i+1}^{n+1} - x_i^{n+1} }{x_{i+1} -  x_i}
	         =  1 + \Delta t \frac{ u_{i+1} - u_i }{x_{i+1} -  x_i}.
\eea
The velocity $u_i$ depends on $u_{c-1}^-, u_{c-1}^+$ and $u_{i+1}$ depends on
$u_{c}^-, u_{c}^+$.
Moreover $u_{c-1}^+ \equiv u_c^-$ by definition. \\
The second fundamental assumption states that the interface normals $\vec{n}_c^\pm$
do not change their direction during their Lagrangian evolution.
The goal is to locate the partial volume into the Lagrangian cell at $t^{n+1}$ and construct the linear interface, knowing its normal $\vec{n}_c^\pm$. Necessarily this partial volume is either in contact with
cell interface $x_i$ (superscript $+$) or $x_{i+1}$ (superscript $-$).
Its volume at $t^{n+1}$ is given by
\bea
V_c^{\pm,n+1} = V_c^\pm \ \dVol_c^\mp
              = V_c^\pm + \Delta t \ \lambda_i^\mp  (u_c^+ - u_c^-).
\eea
If $V_c^{\pm,n+1} \leq V_{i+1/2}^{n+1}$ then there exists a unique line oriented
by the normal  $\vec{n}_c^\pm$ and separating the cell volume into
two sub-volumes $V_c^{\pm,n+1}$ and $(V_{i+1/2}^{n+1}-V_c^{\pm,n+1})$ respectively by the
PLIC (``\textit{Piecewise Linear Interface Construction}'' \cite{DLY1}) method.
As the displacement velocity $u(x)$ is supposed to be piecewise
linear (by the first assumption see Figure~\ref{fig:vitesse}), then, if $x_i < x_c^- < x_{i+1}$ one deduces
$x_i^{n+1} < x_c^{-,n+1} < x_{i+1}^{n+1}$. Therefore  the sub-volume at $t^{n+1}$
is strictly included into the Lagrangian cell volume $V_{i+1/2}^{n+1}$.
This phase is depicted in Figure~\ref{fig:nip-2}-(E)

\subsubsection{{\em Projection} step}
The projection step performs the exact intersection between the Lagrangian condensate
obtained after the reconstruction step in Figure~\ref{fig:nip-2}-(E) and the
Eulerian mesh (bold line squares in Figure~\ref{fig:nip-2}-(A)).
This step is depicted in Figure~\ref{fig:nip-2}-(F). The exact intersection consists
in projecting each partial volume that is accurately located into the condensate, onto some Eulerian fixed cell(s). As instance in
Figure~\ref{fig:nip-2}-(F) the first partial volume is projected onto Eulerian
cells $2$ (green cell) and $3$ (red cell). Contrarily the last partial volume is
totally projected into Eulerian cell $5$ (brown cell).
This projection provides the quantity of material per Eulerian cell, or, equivalently
its volume fraction.

Once volume fractions in the mixed cells are updated though the evolution of
condensates, 2D normals are computed using the same technique as in original
NIP method.

\section{Numerical results} \label{sec:numerics}

In this Section we present a set of test cases to assess the
efficiency of the approach described in the previous Sections.
First, one validates the technique on pure advection test cases
that often present excessive smearing of interfaces due to the
numerical inaccuracy embedded into the scheme.
A square shaped object is advected with constant velocity
 in a diagonal direction in a first test, then into a rotating flow.
Finally an hydrodynamics test case is presented.

\subsection{Advection context} \label{ssec:advection}
\begin{figure}
\begin{center}
\hspace{-1cm}
	\epsfxsize=0.28\textwidth
       \epsfbox{adv_t0}
\hspace{-1cm}
	\epsfxsize=0.28\textwidth
        \epsfbox{adv60_o1}
        \\
\hspace{-1cm}
	\epsfxsize=0.28\textwidth
        \epsfbox{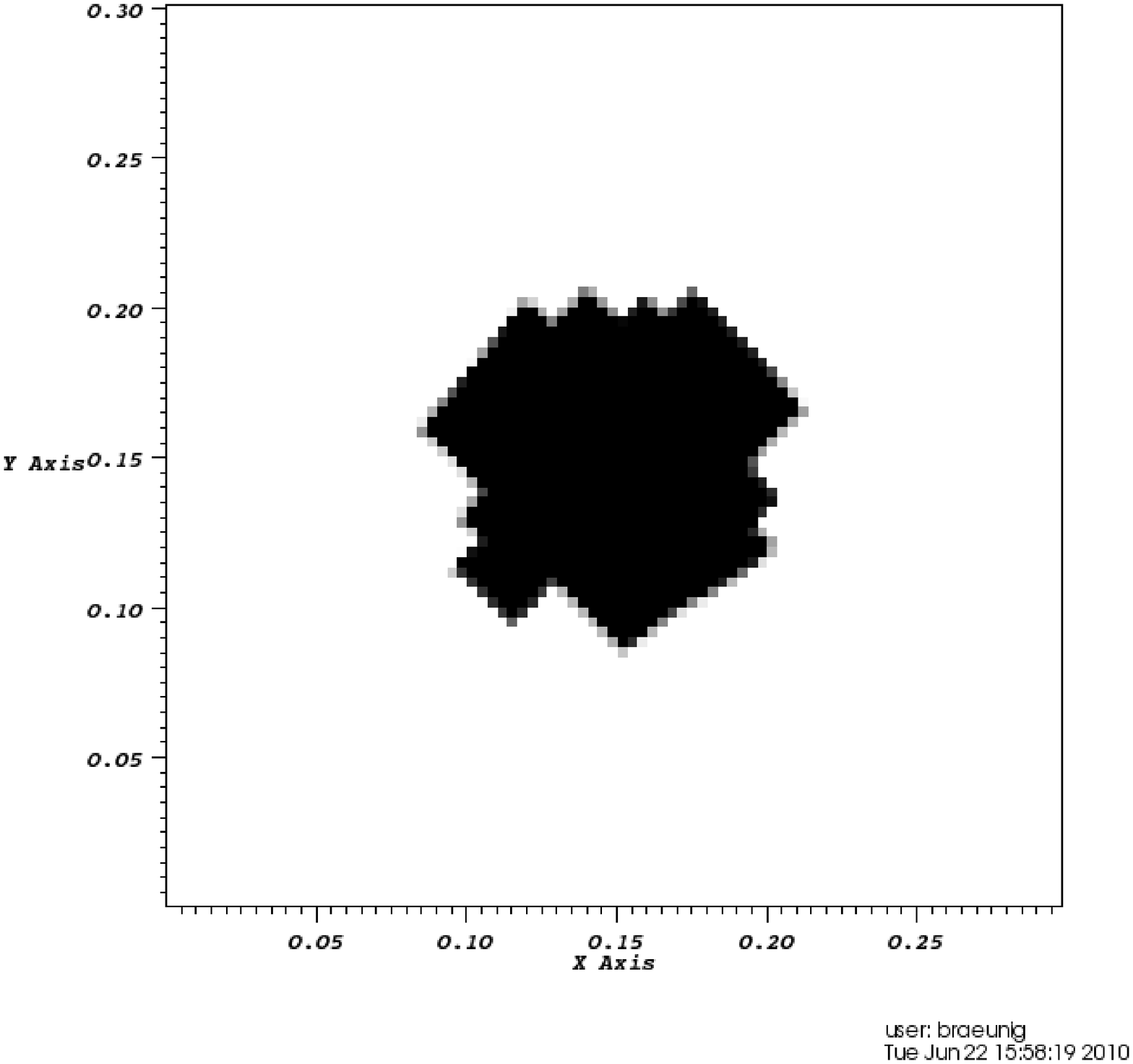}
\hspace{-1cm}
	\epsfxsize=0.28\textwidth
       \epsfbox{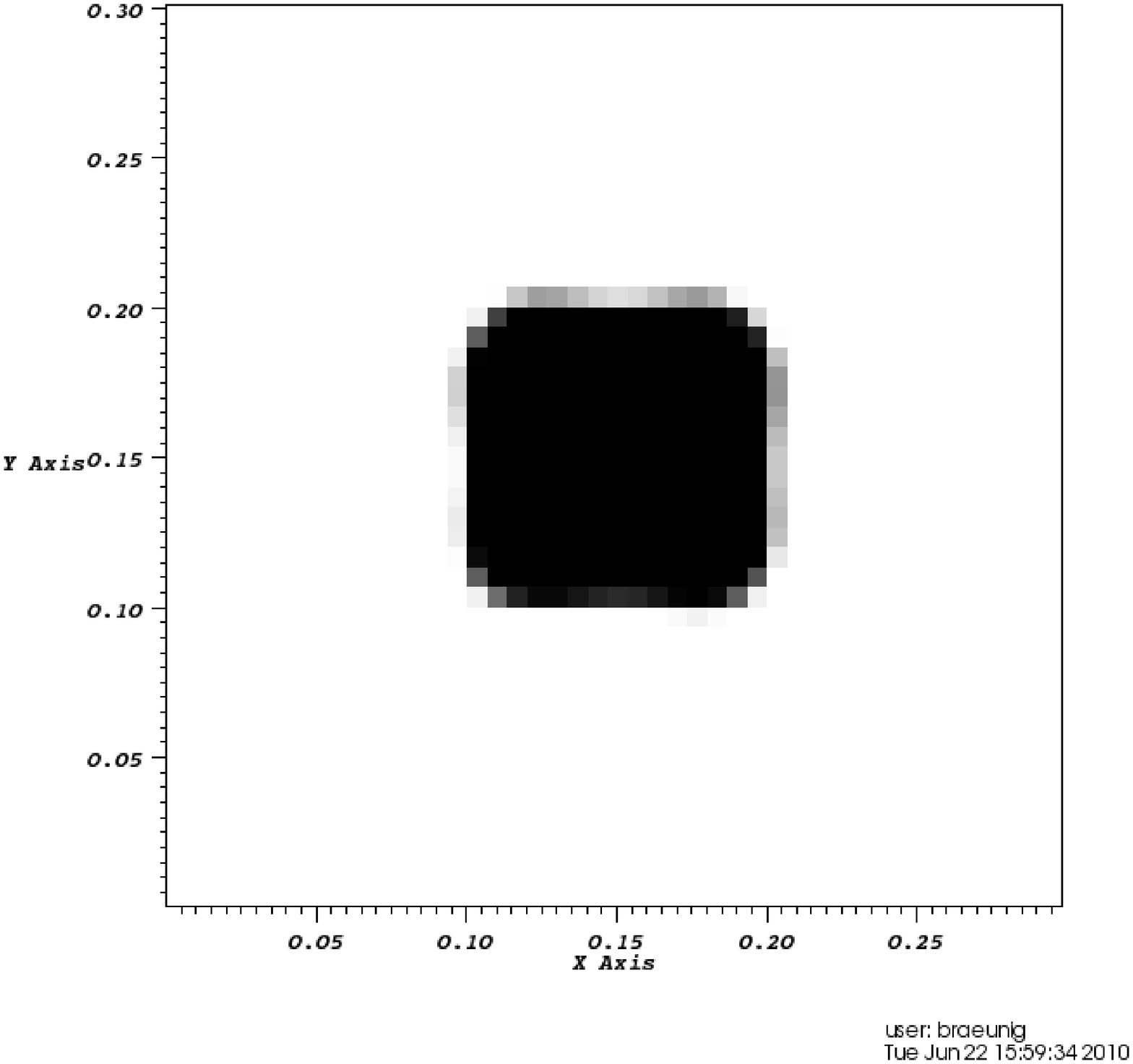}
\caption{ \label{fig:adv_square}
Advection of a square (zoom around the exact position of the initial and final square) ---
Top-left: exact solution ---
Top-right: classical NIP  with a $60\times 60$ mesh ---
Bottom-left: classical NIP  with a $120\times 120$ mesh ---
Bottom-right: ENIP  with a $60\times 60$ mesh.}
\end{center}
\end{figure}

\begin{figure}
\begin{center}
	\epsfxsize=0.5\textwidth
        \epsfbox{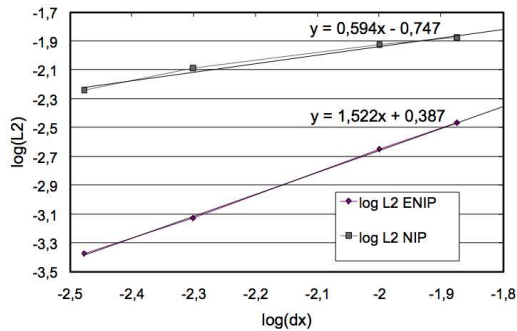}
\caption{ \label{fig:convergence}
Convergence of ENIP vs NIP for a pure advection problem.
The log of the $L_2$ error is displayed as a function of the log of $\Delta x$.}
\end{center}
\end{figure}

An initial square $[0.1;0.1]\times [0.2;0.2]$
is located into the domain $\Omega= [0:0.4]\times [0;0.6]$.
The density into the square is set
to $\rho_0(x)=1$ whereas it is set to $\rho_0(x)=0$ outside.
In the pure advection context this square shape
 should be perfectly conserved through the equation
\bea
\ddt  \rho + u\ddx \rho  + v \ddx \rho = 0,
\eea
where $(u,v)$ is a constant velocity field. The exact solution
at any point $x$ and any time $t$ is
$\rho^{ex}(x,y,t) = \rho_0(x - u \ t, y- v\ t)$. If the numerical
method provides an approximated solution called
$\rho_i^n$ in cell $i$ at time $t^n$ then the error
in $L_\alpha$ norm is evaluated by ($\alpha=1,2$)
\bea
\varepsilon_\alpha = \frac{\sum_i | \rho^n_i - \rho^{ex}(x_i,t^n) |^\alpha}{\sum_i | \rho^{ex}(x_i,t^n) |^\alpha}.
\eea
The first test consists in advecting the square with the constant velocity field $u=1$, $v=3$ up to the time $t=0.1$ then reversing the
advection field by setting $u=-1$, $v=-3$ up to final time $t=0.2$ so that the final configuration perfectly fits the initial one.
Any method (NIP and ENIP included) introduces some error that we intend to measure with this test.
In Figure~\ref{fig:adv_square} are shown the exact solution (top-left) and the results obtained with a $60\times 60$ mesh for NIP (top-right) and ENIP (bottom-right).
ENIP is visibly able to preserve the shape of the square whereas NIP is not. A mesh refinement of NIP computation ($120\times 120$ mesh for the bottom-left panel) does not improve the situation.
In table~\ref{tab:adv} we gather the errors for the $L_1$, $L_2$ norms for successively refined meshes for the
NIP and the proposed ENIP method on this advection problem. Systematically ENIP over-tops NIP.
\begin{table}
\begin{tabular}{|c||cc|cc|}
\hline
$\Delta x=\Delta y$ & \textbf{$L_1$ NIP} & \textbf{$L_1$ ENIP} &  \textbf{$L_2$ NIP} & \textbf{$L_2$ ENIP} \\
\hline
\hline
0.02    &  3.652 & 0.196 & 2.575 & 0.079 \\  
0.0133 &  0.389 & 0.165 & 0.318 & 0.081  \\  
0.01     & 0.339 & 0.111 & 0.284 & 0.053  \\  
0.005    &0.221 & 0.042 & 0.195 & 0.017  \\  
0.0033  &0.155 & 0.025 & 0.138 & 0.010  \\  
\hline
\end{tabular}
\caption{ \label{tab:adv} Error in $L_1, L_2$ norms for the advection problem --- NIP versus ENIP methods.}
\end{table}
In Figure~\ref{fig:convergence} we display the log-log scale results for the error in $L_2$ norm for both methods showing
the improvement gained by ENIP; indeed the slope which represents a measure of the numerical order of convergence is improved by a factor $2.5$ ($0.6$ for NIP and $1.5$ for ENIP).

\begin{figure}
\begin{center}
	\epsfxsize=0.5\textwidth
        \epsfbox{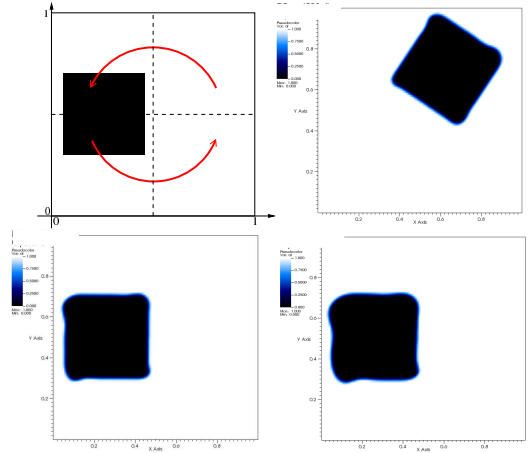}
\caption{ \label{fig:rotation}
Top-left: Sketch of the rigid rotation of a square ---
Top-right: After $5/8$ of the full rotation ---
Bottom-left: After one full rotation ---
Bottom-right: After three full rotations.
}
\end{center}
\end{figure}

The next test consists in the rigid rotation of a square $[0.06;0.46]\times [0.3;0.7]$ (density $1$) into the unit square domain, see
 Figure~\ref{fig:rotation} top-left panel.
A $100\times 100$ uniform mesh is considered and the rotation is given by the velocity field
\bea \nn
u = -100(y-0.5), & & v=100(x-0.5).
\eea
In Figure~\ref{fig:rotation} we display the density  after $5/8$ of the full rotation, after one and three rotations.
The square shape is almost preserved. Contrarily the classical NIP method would totally lose the shape after one rotation.

\subsection{Hydrodynamics context} \label{ssec:physics}
\begin{figure}
\begin{center}
\hspace{-1cm}
	\epsfysize=0.25\textwidth
        \epsfbox{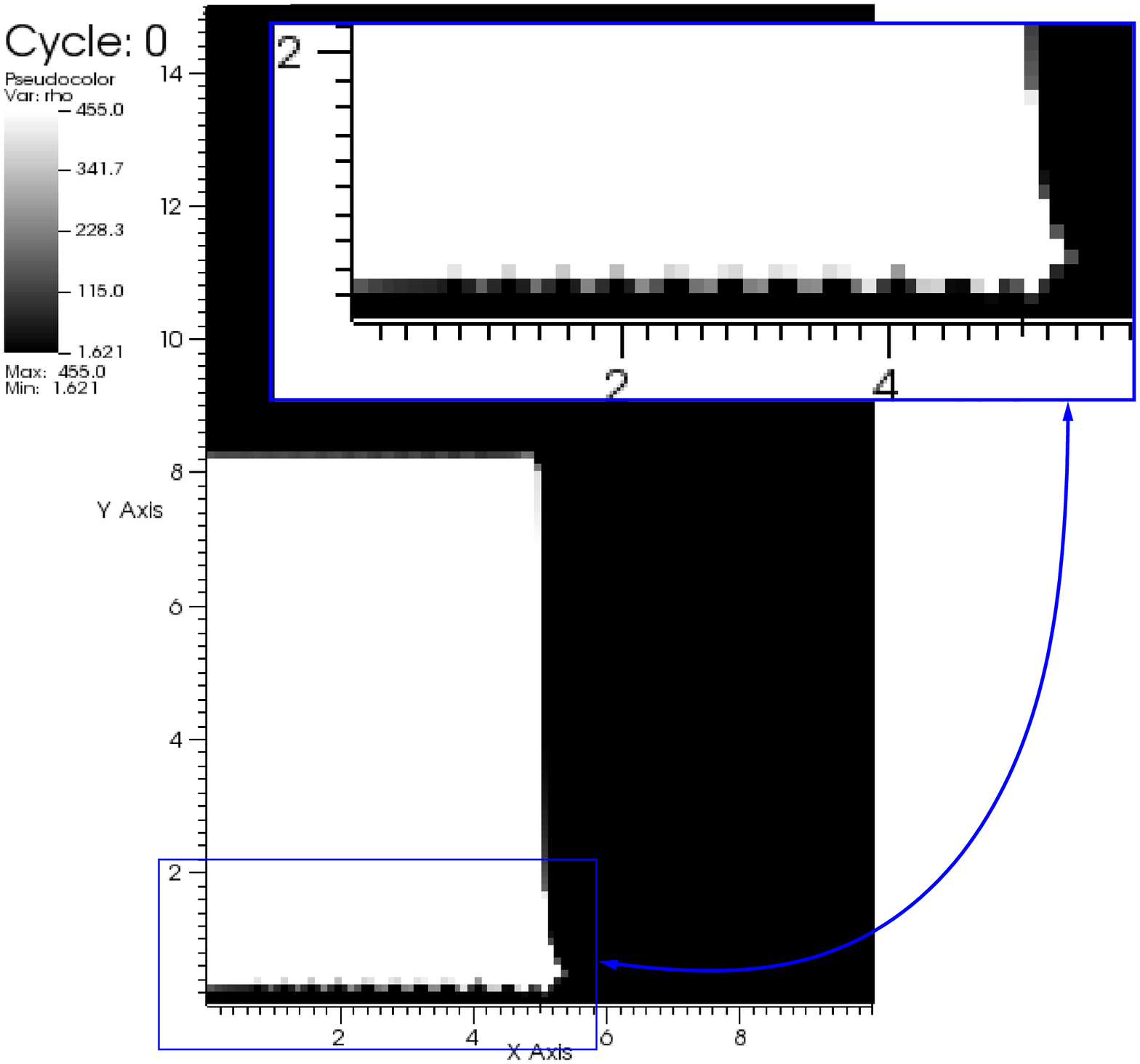}
\hspace{-1cm}
	\epsfysize=0.25\textwidth
        \epsfbox{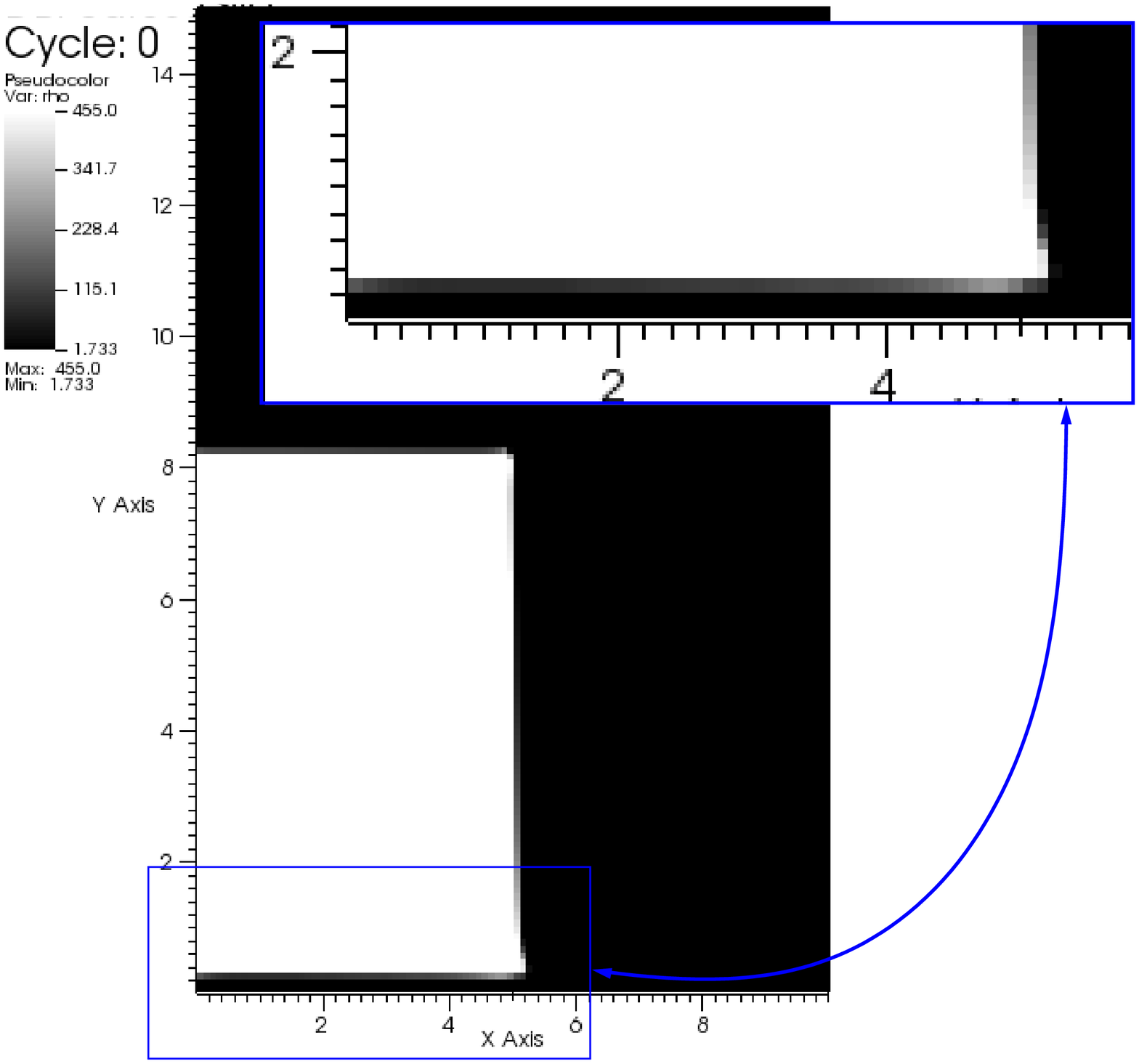}
\\
\textbf{(a)} \hspace{3cm} \textbf{(b)} \\
\hspace{-1cm}
	\epsfysize=0.25\textwidth
        \epsfbox{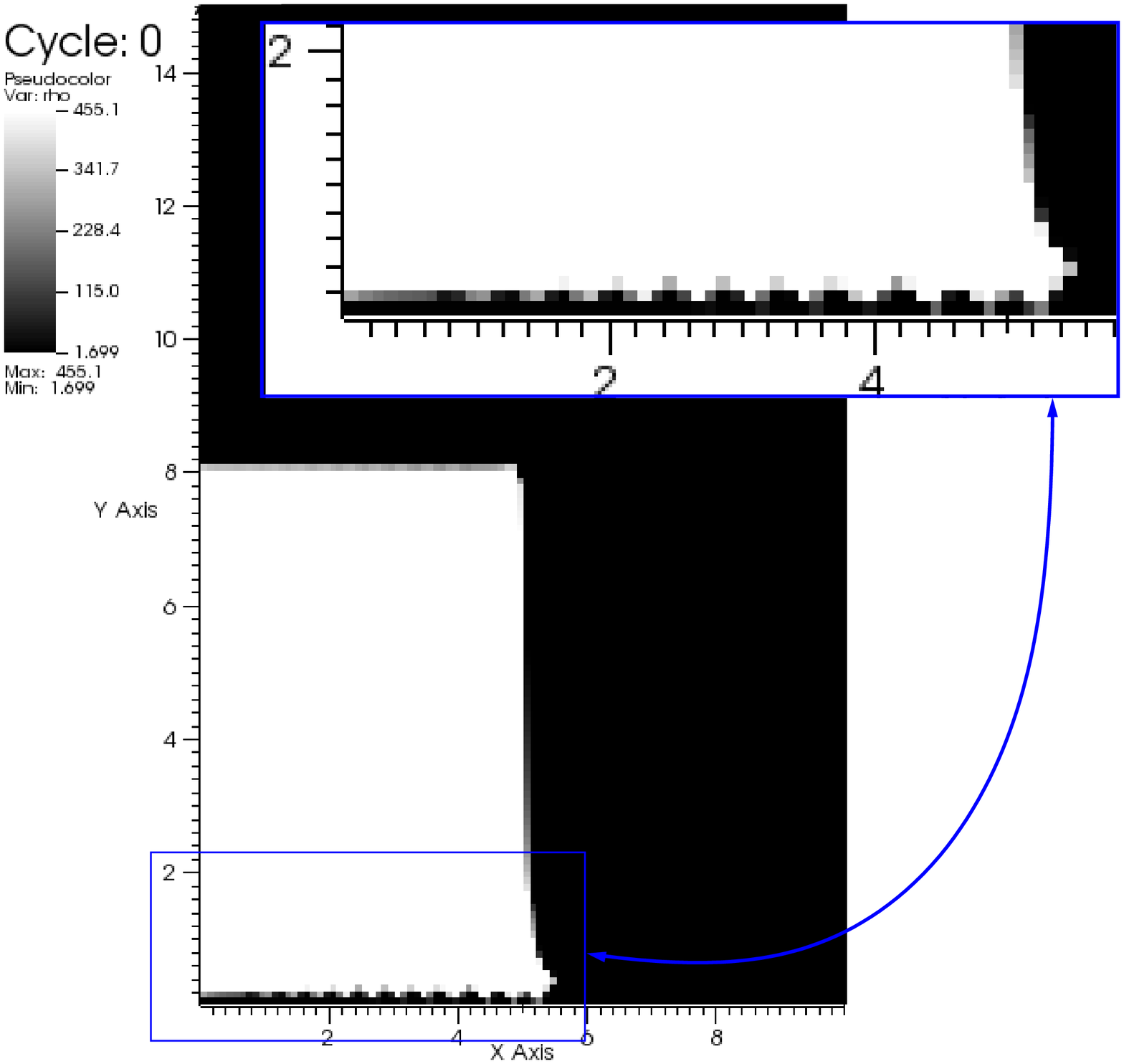}
\hspace{-1cm}
	\epsfysize=0.25\textwidth
        \epsfbox{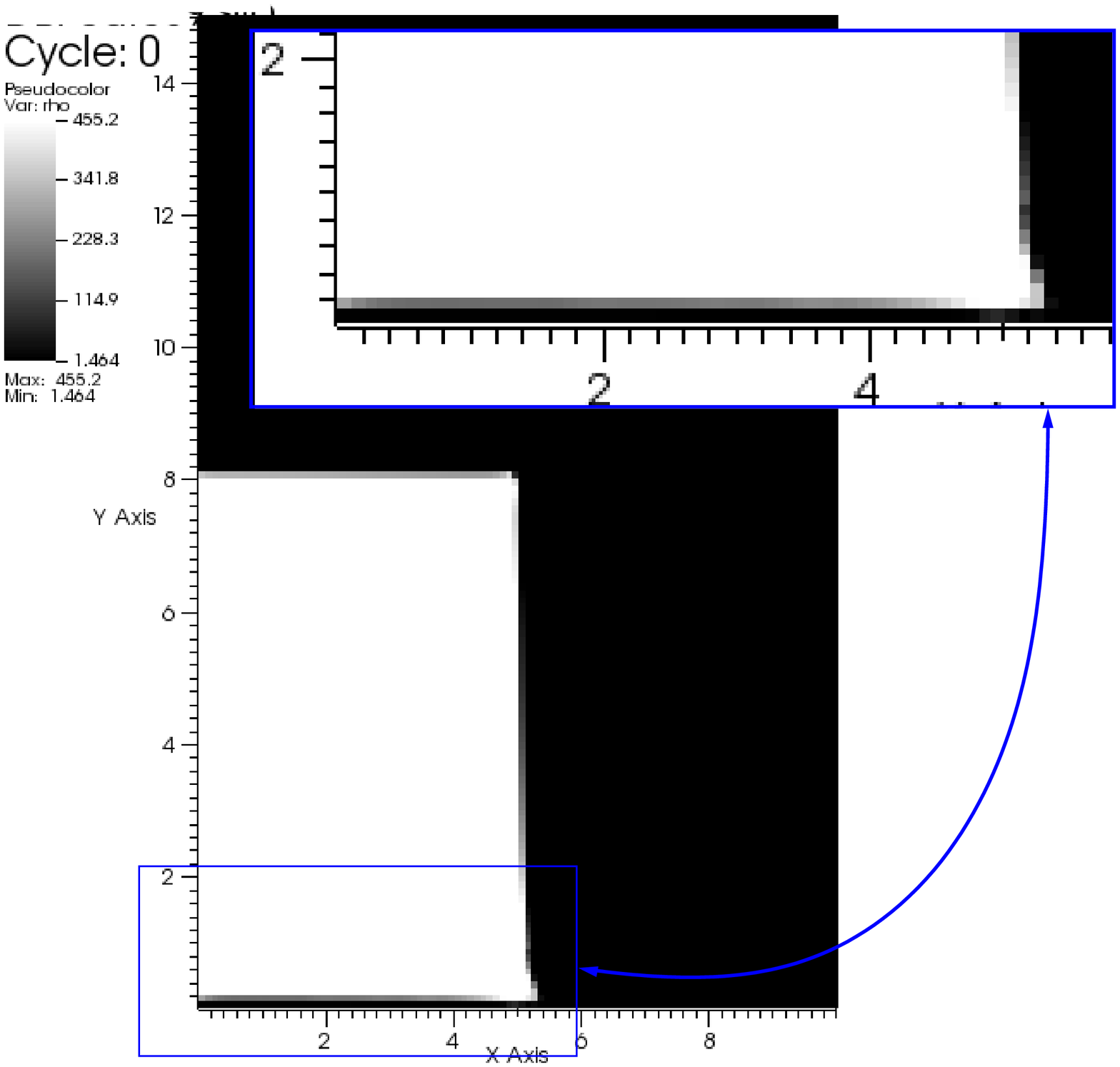}
\\
\textbf{(c)} \hspace{3cm} \textbf{(d)} \\
\caption{ \label{fig:liquid}
---
\textbf{(a)} NIP results at $t=0.6$ (full view and zoom on the impact zone)
---
\textbf{(b)} ENIP results at $t=0.6$
---
\textbf{(c)} NIP results at $t=0.64$
---
\textbf{(d)} ENIP results at $t=0.64$
.}
\end{center}
\end{figure}

We run an idealized 2D test case that corresponds to the free drop of a liquid rectangle within a 2D rectangular tank filled with gas \cite{benchmark}.
This context is inspired by the problem of sloshing that may appear in the tanks of Liquid-Natural-Gas (LNG) carriers.
The study focuses on the ability for the numerical simulations to take properly into account the physics that is of major importance during the liquid impact such as the escape of the gas underneath and its compression.
As a strong sliding process occurs between the compressed gas and the falling liquid. The ability of the method to properly deal with sliding conditions at the interface has a major effect on the final numerical compression and shape of the trapped air. This has ultimately a strong influence on the impact pressure. \\
The test case consists in a domain $\Omega = [0.0;0.0]\times[10m;15m]$ filled with air. The liquid is initially at rest in the rectangle $[0;2]\times[5;10]$ and is falling under the gravity that is pointing downward with magnitude $g=9.81m.s^{-2}$.
A free fall of the liquid into vacuum would impact at $t_{\textrm{impact}}=0.64s$ however due to the presence of the gas this theoretical value is not correct for our simulation however some critical phenomena still occur in the vicinity of this time.  As instance around  $t_{\textrm{impact}}$ a pocket of gas is trapped under the falling liquid and this strongly impacts the numerical impact pressure by decelerating and damping the free fall of the liquid. Therefore a good interface reconstruction method should qualitatively improve the numerical results.
One considers a mesh made of $100\times 150$ uniform cells on the domain.
One shows the results for NIP and ENIP at time $t=0.6s$ Figure~\ref{fig:liquid}-\textbf{(a)-(b)} and $t=0.64s$ in Figure~\ref{fig:liquid}-\textbf{(c)-(d)}.
The classical NIP method was already able to deal with such sliding effects.
However the interface reconstruction method employed is not accurate and stable enough to be free of oscillation that one suspects to be only a numerical artifacts (see panels \textbf{(a-c)}).
Contrarily the new reconstruction method ENIP on this very same test case is able to produce a smooth interface that permits to obtain a more realistic simulation.  Indeed this simulation prominently displays the fact that the ``bubbling'' effects of NIP is of pure numerical origin and that ENIP cures this drawback.

\section{Conclusion and perspectives} \label{sec:conclusion}

This paper deals with the improvement of the so-called NIP (Natural Interface Positioning) method.
The NIP method described in \cite{Braeunig09} is an add-on to the FVCF method in order to treat multi-material fluid flows
uses the concept of condensate. A condensate is the association of contiguous mixed cells in either $x$ or $y$ direction. They are further treated as an entity to make possible the treatment of each mixed cell taken individually. NIP is the method based on the following steps:
 {\em Representation}, {\em Condensate construction},  {\em Condensate evolution}, {\em Reconstruction},  and {\em Projection}.
The present paper points the weakness of the NIP method in pure advection context and, consequently, in a full multi-material hydrodynamics one. An enhanced NIP method is proposed (ENIP). It modifies several of the previous listed steps. More precisely the condensate is assumed to evolved in an almost-Lagrangian fashion. The reconstruction step assumes that the condensate keeps the same form modulo some expansion/compression that the numerical scheme already provides. So the displacement of the condensate is performed either with the true computed velocity or with an interpolation of it. \textit{In fine} the condensate preserves its topology contrarily to the original NIP method for which the condensate has no recollection of its shape from the beginning of the time step. \\
The capability of the full numerical method is now dramatically improved as seen on advection test cases (advection and rigid rotation of a square). Moreover we ran ENIP on a difficult mutli-material hydrodynamics tests simulating the free drop of a liquid rectangle within a 2D rectangular tank filled with gas
in the context of sloshing that may appear in the tanks of Liquid-Natural-Gas carrier (see \cite{benchmark}). The accuracy, stability and robustness of the ENIP method is clearly seen especially at the time some air is trapped under the water.
In the near future we plan to investigate the evolution of this method to the case of mixed cells with more than two materials. In this case the only difficulty lays in the positioning of the different materials in the cell, but their evolution within the condensate follows exactly the same algorithm ENIP with no modification of the numerical scheme. We also plan to investigate the evolution of the method in 3D.

 

\end{document}